\documentclass[preprint,12pt]{elsarticle}

\usepackage{amsmath,amssymb,amsthm}
\usepackage{algorithm}
\usepackage{algpseudocode}
\usepackage{graphicx}
\usepackage{subcaption}
\usepackage{tikz}
\usepackage{float}
\usepackage{stmaryrd}        % \llbracket, \rrbracket
\usepackage{xcolor}
\usepackage[colorlinks=true,
            linkcolor=blue,
            citecolor=blue,
            urlcolor=blue]{hyperref}

\newtheorem{remark}{Remark}[section]

\newcommand{\fatu}{\mathbf{u}}
\newcommand{\fatw}{\mathbf{w}}
\newcommand{\fatn}{\mathbf{n}}
\newcommand{\Th}{{\mathcal{T}_{\!h}}}
\newcommand{\Eh}{\mathcal{E}_h}
\newcommand{\dx}{\,\mathrm{d}x}
\newcommand{\dS}{\,\mathrm{d}S}
\newcommand{\drm}{\mathrm{d}}
\newcommand{\Hdiv}{H(\mathrm{div})}

\journal{Journal of Computational Physics}

\begin{document}

\begin{frontmatter}

\title{A Second-Order Semi-Implicit Projection Time-Stepping Scheme for the
       Shallow-Water Equations on the Sphere Using Compatible
       Finite Elements}

%% -- Authors (fill in) --------------------------------------
\author[1,2]{Karina Kowalczyk\corref{cor}}
\author[1]{Colin Cotter}

\cortext[cor]{Corresponding author}
\address[1]{Department of Mathematics, Imperial College London, UK}
\address[2]{Potsdam Institute for Climate Impact Research, Germany}

%% -- Abstract -----------------------------------------------
\begin{abstract}
We propose a second-order, semi-implicit time-stepping scheme for the shallow-water equations on the sphere discretised with compatible finite elements. The scheme splits the advection operator from the wave operator, applying a DG-upwind SSPRK3 method to the advection part and an $\Hdiv$-conforming implicit discretisation to the wave part, solved efficiently via hybridisation. The scheme follows a three-step symmetric splitting: an explicit half-step, a full advection step, and an implicit half-step that simultaneously projects the velocity back into the $\Hdiv$ space. This symmetry guarantees global second-order accuracy and confines the stability restriction to the CFL condition of the advection step alone. The time-step can therefore be extended beyond this limit by subcycling the advection step without additional implicit solves. The method is implemented in Firedrake and validated against a standard suite of spherical shallow-water benchmarks, demonstrating stable results for Courant numbers up to 1 (with two advection subcycles) with second-order finite elements, second-order convergence, and good parallel scaling.
\end{abstract}


\end{frontmatter}

%% ============================================================
%%  Sections
%% ============================================================
%% ============================================================
%%  Section 1 — Introduction
%% ============================================================
\section{Introduction}
\label{sec:intro}

Numerical weather prediction (NWP) and climate projection depend
critically on efficient and accurate discretisations of the
equations governing atmospheric and oceanic flow.  A central
challenge is the multiscale nature of these flows: fast
acoustic and gravity waves coexist with slow Rossby waves and
nonlinear advective dynamics, spanning many orders of magnitude
in wave speed.  The design of time integration schemes that
handle this efficiently, while remaining accurate,
stable, and scalable on massively parallel hardware, is an
active area of research.

Semi-implicit methods are the dominant paradigm in operational
NWP: the fast wave modes are treated implicitly, removing the
associated upper bound  on the time step, while the slower dynamics are
advanced explicitly or semi-explicitly.  Classical semi-implicit
approaches date to the work of Robert, Henderson and
Turnbull \cite{robert1972implicit} and have since been refined
extensively. In operational forecasting, the semi-implicit
treatment is commonly combined with a semi-Lagrangian advection
scheme, yielding the semi-implicit semi-Lagrangian (SISL)
approach that allows time steps far exceeding the advective
Courant-Friedrichs-Levi (CFL) limit \cite{melvin2010inherently}.  Alternative splitting
strategies include horizontally explicit, vertically implicit
(HEVI) methods \cite{weller2013runge,lock2014numerical} and
split-explicit schemes that use subcycled timestepping for the fast
modes \cite{klemp1978numerical,skamarock1992efficiency}, each
offering different trade offs between stability, accuracy,
and parallel scalability.

Compatible finite element methods, underpinned by the finite element exterior
calculus \cite{arnold2006finite, arnold2018finite,arnold2010finite}, have
gained significant traction as a foundation for next generation
atmospheric dynamical cores.  The compatible finite element
framework, proposed in the context of geophysical fluid dynamics by Cotter and
Shipton \cite{cotter2012mixed} and further developed in
\cite{cotter2014finite,natale2016compatible}, uses a hierarchy
of finite element spaces that form a discrete de Rham
complex.  In particular, the velocity lies in an
$\Hdiv$-conforming space and the depth (or pressure) in a
discontinuous $L^2$ space, mirroring the continuous
structure of the equations and enabling discrete analogues of
key vector calculus identities.  This framework underpins
the GungHo dynamical core being developed for the UK Met
Office's next generation NWP
model \cite{melvin2019mixed,melvin2024mixed}, as well as
research codes such as Gusto \cite{gusto2023} built on the
Firedrake finite element system \cite{FiredrakeUserManual}.

In the Met Office's semi-implicit
discretisation \cite{melvin2019mixed}, advection is handled via
a finite volume scheme, requiring an explicit mapping between the
finite volume and finite element representations at each
time step.  Within each time step the Newton linearisation is
split into outer and inner loops, with the advection scheme
invoked in the outer loop and the implicit wave solve in both,
resulting in multiple expensive residual evaluations per
time step.  Formulating a purely finite element time stepping
scheme that reduces these evaluations, maintains second order
accuracy, and admits scalable solvers is therefore of direct
practical interest.

The theoretical foundation of the present work is the
projection method framework developed by Klein, Vater, Benacchio
and collaborators for finite volume
models \cite{vater2009stability,vater2018semi,benacchio2019semi,
benacchio2014blended,klein2016doubly,chew2022one}.
In \cite{vater2009stability}, Vater and Klein established a
projection method for zero Froude number (incompressible limit)
shallow water flows, combining an explicit predictor step for an
auxiliary hyperbolic system with two divergence projection
correction steps.  This was extended to low Froude number flows
via a semi-implicit formulation in \cite{vater2018semi} and
extended to  compressible finite volume models by Benacchio
and Klein in \cite{benacchio2019semi}, yielding a symmetric,
three step splitting that wraps an explicit advection step
between explicit and implicit Euler half steps of the wave part.
The resulting scheme is second order accurate and its stability
is governed solely by the CFL condition of the advection step.
Hastermann formulated the projection step of these schemes by
considering mixed finite elements, providing a step to the present
work \cite{hastermann2022analysis}.

We translate this three step framework to the compatible finite
element setting for the shallow water equations on the sphere.
The shallow water equations (SWEs) constitute the simplest
model exhibiting the slow-fast split of atmospheric dynamics
and are a standard first tier of the model hierarchy used to
validate dynamical core components \cite{williamson1992standard}.
Splitting the advection from the wave part allows us to apply a
Discontinuous Galerkin (DG) upwind finite element discretisation \cite{cockburn2001runge} for the
advection, which is well suited to the discontinuous spaces and
is known to suppress spurious oscillations through its
dissipative upwind flux.  The wave part is handled with an
$\Hdiv$-conforming discretisation and solved efficiently via
hybridisation \cite{cockburn2009unified}, which reduces the
global implicit system to a smaller, sparse condensed problem on
element facets.  The implicit final step simultaneously advances
the wave dynamics and projects the discontinuous advected
velocity back into the $\Hdiv$ space, in direct analogy to the
projection steps in \cite{vater2009stability,vater2018semi}.

The resulting scheme is a purely Eulerian, fully finite element
method requiring no finite volume to finite element mappings
and no nested solver loops for stabilisation. The advecting
velocity is approximated to first order at the half time step
via a two step subscheme. Following \cite{benacchio2019semi},
this suffices to maintain overall second order accuracy. The stability of the full scheme
is restricted only by the advection CFL condition; by subcycling
the advection step $N$ times with time step $\Delta t/N$, the
effective Courant number of the advection step is reduced by a factor of $N$, allowing larger overall time steps without additional implicit solves.

The paper is organised as follows.
Section \ref{sec:swe_sphere} presents the governing equations
and their advection-wave split.
Section \ref{sec:method} introduces the three step time stepping
scheme in strong form, establishes second order accuracy,
describes the spatial finite element discretisation including
the DG bilinear forms, and details scalable solver strategies
for each step.
Section \ref{sec:results} presents numerical results for three
standard benchmark test cases.
Section \ref{sec:conclusion} concludes with a summary of
findings and directions for future work.

%% ============================================================
%%  Section 2 — Governing Equations
%% ============================================================
\section{Governing Equations}
\label{sec:swe_sphere}

The shallow water equations (SWEs) arise from the vertical
integration of the three dimensional Navier-Stokes equations
under the assumption that the characteristic horizontal length
scale significantly exceeds the fluid depth, yielding a
two-dimensional model that nevertheless captures the essential
physics of large scale geophysical flows, including gravity
waves, Rossby waves, and geostrophic adjustment.  We consider
the SWEs on the sphere, where the flow is constrained to the
surface of a sphere of radius $R$.  In a rotating frame,
including bottom topography $b$, the equations read
\begin{align}
    \frac{\partial \mathbf{u}}{\partial t}
        + (\fatu \cdot\nabla)\fatu
        + f\mathbf{u}^{\perp}
        &= -g\nabla (D+b),
    \label{eq:SWE_u}\\[4pt]
    \frac{\partial D}{\partial t}
        + \nabla \cdot (\mathbf{u}D)
        &= 0,
    \label{eq:SWE_D}
\end{align}
where $\fatu$ is the horizontal velocity field tangent to the
sphere, $D$ is the fluid column depth, $g$ is the gravitational
acceleration, and
\begin{equation}
    f = \frac{2\Omega z}{R}
\end{equation}
is the Coriolis parameter for Earth's angular velocity $\Omega$
and position vector $\mathbf{r} = (x,y,z)$.  The vector
\begin{equation}
    \fatu^\perp = (\mathbf{r} \times \fatu)/R
\end{equation}
is the velocity perpendicular to $\fatu$ in the tangent plane,
and $\nabla$ denotes the projection of the three dimensional
gradient into the tangent plane. 
%When solved on the surface of a sphere, the velocity equation \eqref{eq:SWE_u} is understood to be projected onto the tangent plane, so that u remains tangent to the sphere for all time.

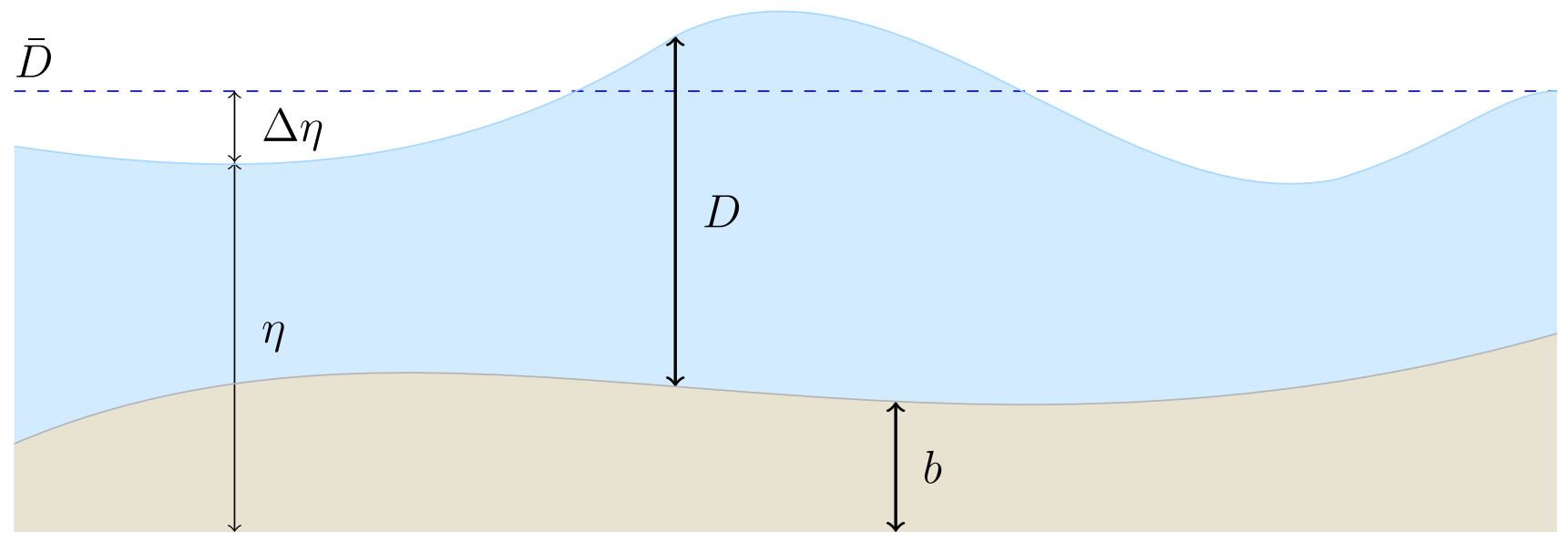
\begin{figure}[H]
\begin{center}
\scalebox{0.8}{
\begin{tikzpicture}
    \definecolor{waterblue}{RGB}{165,216,255}
    \definecolor{sandbeige}{RGB}{232,226,208}
    \definecolor{lightgray}{RGB}{180,180,180}
    % Water body
    \fill[waterblue!50]
        (0,3.5) .. controls (2,3.2) and (4,3.2) .. (6,4.5)
                .. controls (8,5.5) and (10,2.8) .. (12,3.2)
                .. controls (13,3.5) and (13.5,4) .. (14,4)
                -- (14,0) -- (0,0) -- cycle;
    % Bottom topography
    \fill[sandbeige]
        (0,0.8) .. controls (4,2.5) and (8,0.1) .. (14,1.8)
                -- (14,0) -- (0,0) -- cycle;
    % Mean water level
    \draw[blue, dashed] (0,4) -- (14,4);
    \node[left] at (0.5,4.3) {$\Bar{D}$};
    % Depth D
    \draw[thick,<->] (6,4.5) -- (6,1.32);
    \node[right] at (6.1,2.91) {$D$};
    % Free surface perturbation
    \draw[<->] (2,3.355) -- (2,4);
    \node[right] at (2.1,3.65) {$\Delta\eta$};
    \draw[<->] (2,3.355) -- (2,0);
    \node[right] at (2.1, 1.77) {$\eta$};
    % Topography height b
    \draw[thick,<->] (8,1.18) -- (8,0);
    \node[right] at (8.1,0.59) {$b$};
    % Water surface curve
    \draw[waterblue]
        (0,3.5) .. controls (2,3.2) and (4,3.2) .. (6,4.5)
                .. controls (8,5.5) and (10,2.8) .. (12,3.2)
                .. controls (13,3.5) and (13.5,4) .. (14,4);
    % Topography curve
    \draw[lightgray]
        (0,0.8) .. controls (4,2.5) and (8,0.1) .. (14,1.8);
\end{tikzpicture}
}
\end{center}
\caption{Schematic of a shallow water system with reference depth
$\Bar{D}$, free surface $\eta$, fluid column depth $D$, and
bottom topography $b$.}
\label{fig:SWE}
\end{figure}

A schematic of the shallow water system is shown in Figure \ref{fig:SWE}.

\subsection{Advection-wave splitting}
\label{subsec:splitting}

To separate the advective dynamics from the wave dynamics, let
$\Bar{D}$ be a reference depth (in practice taken as the
constant mean depth; see Remark \ref{rem:Dbar} below).
Adding and subtracting $\nabla\cdot(\fatu\Bar{D})$ in the
continuity equation, the system \eqref{eq:SWE_u}-\eqref{eq:SWE_D}
can be rewritten as
\begin{align}
    \frac{\partial \mathbf{u}}{\partial t}
        + (\fatu \cdot\nabla)\fatu
        + f\mathbf{u}^{\perp}
        + g\nabla (D+b)
        &= 0,
    \label{eq:SWE_split_u}\\[4pt]
    \frac{\partial D}{\partial t}
        + \nabla \cdot \!\big(\fatu(D-\Bar{D})\big)
        + \nabla \cdot (\mathbf{u}\bar{D})
        &= 0.
    \label{eq:SWE_split_D}
\end{align}
This decomposition isolates an advective contribution
$\nabla\cdot(\fatu(D-\Bar{D}))$, which transports the depth
perturbation with the flow, from a wave like contribution
$\nabla\cdot(\fatu\Bar{D})$, which drives the fast gravity
wave dynamics. The time stepping scheme introduced in
Section \ref{sec:method} exploits this split, advancing each
part with the most appropriate discretisation.

For constant $\Bar{D}$, discarding the advective term
$\nabla\cdot(\fatu(D-\Bar{D}))$ yields the linearised
shallow water equations,
\begin{align}
    \frac{\partial \mathbf{u}}{\partial t}
        + f\mathbf{u}^{\perp}
        + g\nabla (D+b) &= 0,\\[4pt]
    \frac{\partial D}{\partial t}
        + \nabla \cdot (\mathbf{u}\bar{D}) &= 0,
\end{align}
which exhibit purely wave dynamics. The lower bound for the wave phase velocity is
$c = \sqrt{g\Bar{D}}$.

\begin{remark}
\label{rem:Dbar}
Although mathematically $\Bar{D}$ can be any reference profile,
stability issues were observed in our experiments when $\Bar{D}$
was not chosen carefully.  Choosing $\Bar{D}$ as the constant
mean depth is robust across all test cases considered here.
\end{remark}

%% ============================================================
%%  Section 3 — Numerical Method
%% ============================================================
\section{Numerical Method}
\label{sec:method}

\subsection{Time discretisation}
\label{subsec:time-disc}

Given the solution $(\fatu^{n}, D^{n})$ at time $t^n$, we seek
$(\fatu^{n+1}, D^{n+1})$ at $t^{n+1} = t^n + \Delta t$.  The
proposed scheme is a symmetric three step splitting method
inspired by the finite volume framework of
Benacchio and Klein \cite{benacchio2019semi}, consisting 
of the following steps.
\begin{enumerate}
    \item An explicit (forward Euler) half time step for the
          wave part, using an $\Hdiv$-conforming velocity space;
    \item a full time step of DG upwind advection, using a
          three stage SSPRK3 scheme in time and a broken
          $\Hdiv$ space for the velocity; and
    \item an implicit (backward Euler) half time step for the
          wave part, using an $\Hdiv$-conforming velocity space.
\end{enumerate}
The finite element spaces mentioned here are discussed in Section \ref{subsec:space-disc}.

Steps 1 and 3 together form a Crank-Nicolson discretisation
of the wave operator, which is unconditionally
stable \cite{durran2010numerical}.  Because the velocity is
advanced through a discontinuous space in Step 2, the
advected velocity field is in general not $\Hdiv$ conforming.
Step 3 simultaneously advances the wave dynamics and restores
$\Hdiv$ conformity, in direct analogy to the divergence
projection steps in \cite{vater2009stability,vater2018semi};
it is therefore referred to as the \emph{projection step}.

The advecting velocity $\bar{\fatu}$ used in Step 2 is a
first order approximation to $\fatu(t^{n+1/2})$,
the solution at the half time $t^{n+1/2} = t^n + \Delta t/2$.
As shown in Section \ref{subsec:advecting_velocity}, this
is sufficient to maintain overall second order accuracy.
The full algorithm is summarised in
Algorithm \ref{alg:timestepping}.

\medskip
\noindent\textbf{Step 1: Explicit step.}
Let $(\hat{\fatu}, \hat{D})$ be the solution of the
forward Euler half step for the wave part, according to
\begin{gather}
    \hat{\mathbf{u}} - \mathbf{u}^n
    + \frac{\Delta t}{2} f (\mathbf{u}^n)^{\perp}
    = - g\, \frac{\Delta t}{2}\nabla(D^n+b),
    \label{eq:explicit_u}\\[4pt]
    \hat{D} - D^n
    + \frac{\Delta t}{2}\nabla\cdot(\mathbf{u}^n\Bar{D})
    = 0.
    \label{eq:explicit_D}
\end{gather}

\medskip
\noindent\textbf{Step 2: Advection step.}
Let $(\fatu^*, D^*)$ be the output of the SSPRK3 scheme applied
to the advection equations over a full time step $\Delta t$
with frozen advecting velocity $\Bar{\fatu}$, according to
\begin{align}
    \fatu^{(1)} &= \hat{\fatu}
        - \Delta t \, (\Bar{\fatu}\cdot \nabla)\hat{\fatu},\\
    \fatu^{(2)} &= \tfrac{3}{4}\hat{\fatu}
        + \tfrac{1}{4}\!\left(\fatu^{(1)}
        - \Delta t\, (\Bar{\fatu}\cdot \nabla)\fatu^{(1)}\right),\\
    \fatu^* &= \tfrac{1}{3}\hat{\fatu}
        + \tfrac{2}{3}\!\left(\fatu^{(2)}
        - \Delta t\, (\Bar{\fatu}\cdot \nabla)\fatu^{(2)}\right),
\end{align}
and analogously for the depth,
\begin{align}
    D^{(1)} &= \hat{D}
        - \Delta t\, \nabla \cdot \!\left(
          \Bar{\fatu} \, (\hat{D} - \Bar{D})\right),\\
    D^{(2)} &= \tfrac{3}{4}\hat{D}
        + \tfrac{1}{4}\!\left(D^{(1)}
        - \Delta t\,\nabla \cdot \!\left(
          \Bar{\fatu}(D^{(1)} - \Bar{D})\right)\right),\\
    D^* &= \tfrac{1}{3}\hat{D}
        + \tfrac{2}{3}\!\left(D^{(2)}
        - \Delta t\,\nabla \cdot \!\left(
          \Bar{\fatu}(D^{(2)} - \Bar{D})\right)\right).
\end{align}
A three stage Runge-Kutta method is used because the DG
spatial discretisation (Section \ref{subsec:space-disc})
uses quadratic ($p=2$) polynomials: Cockburn and
Shu \cite{cockburn2001runge} showed that SSPRK3 is required
for $L^2$-stability of DG upwind schemes with $p=2$.

\medskip
\noindent\textbf{Step 3: Projection step.}
Let $(\fatu^{n+1}, D^{n+1})$ be the solution of the
backward Euler half step for the wave part, according to
\begin{gather}
    \mathbf{u}^{n+1} - \mathbf{u}^*
    + \frac{\Delta t}{2} f (\mathbf{u}^{n+1})^{\perp}
    = - g\, \frac{\Delta t}{2} \, \nabla (D^{n+1}+b),
    \label{eq:proj_u}\\[4pt]
    D^{n+1} - D^*
    + \frac{\Delta t}{2}\nabla\cdot(\mathbf{u}^{n+1}\Bar{D})
    = 0.
    \label{eq:proj_D}
\end{gather}
This step restores $\Hdiv$ conformity of the velocity,
projecting $\fatu^*$ (which lies in the broken space
$\hat{V}_h^1$) back into the conforming space $V_h^1$.

\begin{algorithm}[H]
\caption{Three step projection time stepping scheme}
\label{alg:timestepping}
\begin{algorithmic}
\Require $(\fatu^n, D^n)$ at time $t^n$; subcycle count
         $N\in\mathbb{N}$; advecting velocity $\bar{\fatu}$
\Ensure  $(\fatu^{n+1}, D^{n+1})$ at time $t^{n+1}=t^n+\Delta t$
\State \textbf{Explicit step}: solve for $(\hat{\fatu},\hat{D})$
    via \eqref{eq:explicit_u}-\eqref{eq:explicit_D}
    \hfill [$H(\mathrm{div})$ space]
\State \textbf{Advection step}: \textbf{for} $k=1,\dots,N$
    \textbf{do}
\State \quad solve SSPRK3 with time step $\delta t = \Delta t/N$
    and advecting velocity $\bar{\fatu}$
    \hfill [broken $H(\mathrm{div})$ space]
\State \quad $(\hat{\fatu},\hat{D})\gets(\fatu^*,D^*)$
\State \textbf{end for}
\State \textbf{Projection step}: solve for $(\fatu^{n+1},D^{n+1})$
    via \eqref{eq:proj_u} -\eqref{eq:proj_D}
    \hfill [$H(\mathrm{div})$ space]
\State \Return $(\fatu^{n+1}, D^{n+1})$
\end{algorithmic}
\end{algorithm}

\subsection{Subcycling the advection step}
\label{subsec:subcycling}

The Courant number for the advection step is
\begin{equation}
    C = \max_K \frac{\Delta t}{|K|}
        \int_{\partial K} \max(-\fatu\cdot\mathbf{n},\,0)\,\mathrm{d}s ,
    \label{eq:courant}
\end{equation}
{which measures the inflow of volume through the boundary of a cell $K$ during
one time step, relative to its area $|K|$. In
one dimension it reduces to the usual $C = |\fatu|\,\Delta t/\Delta x$.}

{The DG-upwind scheme is conditionally stable, requiring
$C \leq C_{\mathrm{crit}}$, for some scheme dependent $C_{\mathrm{crit}}>0$. {For degree $p = 2$ with
SSPRK3 time stepping, the one dimensional $L^2$ stability analysis of
\cite{cockburn2001runge} gives $C_{\mathrm{crit}} = 0.209$.}
}
To address this, since the wave part (Steps 1 and 3 together) {forms a
Crank--Nicolson scheme and} is unconditionally stable, one may take a
larger global time step $\Delta t$ and subcycle the advection step.
Instead of one advection step with time step $\Delta t$, we perform $N$
consecutive steps each with time step $\delta t = \Delta t/N$, reducing
the effective Courant number of the advection step by a factor of $N$ while keeping the overall
time step fixed at $\Delta t$. This allows the use of larger time steps
than the advection CFL condition would otherwise permit, without
increasing the number of implicit projection solves.
%Throughout Section~\ref{sec:results}, the quoted Courant number refers to
%the effective Courant number in the advection sub steps, that is
%{$C_{\mathrm{eff}} = C/N$, i.e.\ \eqref{eq:courant} with
%$\Delta t$ replaced by $\delta t$}.

\subsection{Second order accuracy}
\label{subsec:accuracy}

Although Steps 1 and 3 are individually first order, the
full three step scheme achieves global second order by symmetry.
Let $\Phi_{\Delta t/2}$ denote the forward Euler step,
$\Phi^*_{\Delta t/2}$ its adjoint (backward Euler step), and
$\Psi_{\Delta t}$ the SSPRK3 advection step.  The method reads
\begin{equation}
    \mathcal{M}_{\Delta t}
    = \Phi_{\Delta t/2} \circ \Psi_{\Delta t}
      \circ \Phi^*_{\Delta t/2}.
\end{equation}
A direct comparison of coefficients gives
\begin{equation}
    \Psi_{\Delta t}
    = \Psi_{\Delta t/2} \circ \Psi^*_{\Delta t/2}
      + \mathcal{O}(\Delta t^3),
\end{equation}
and thus, setting $\Gamma := \Phi \circ \Psi$,
\begin{equation}
    \mathcal{M}_{\Delta t}
    = \Gamma_{\Delta t/2} \circ \Gamma^*_{\Delta t/2}
      + \mathcal{O}(\Delta t^3).
\end{equation}
The method is therefore symmetric up to $\mathcal{O}(\Delta
t^3)$ terms and has global order two
by \cite[Theorem 3.2]{hairer2006geometric}.  An alternative
argument \cite{benacchio2019semi} notes that the
Crank--Nicolson rule is effectively applied along the
Lagrangian trajectories implicitly traced by the advection
scheme.

\subsection{The advecting velocity}
\label{subsec:advecting_velocity}

Following \cite{benacchio2019semi} (Appendix A), a first order
accurate approximation to $\fatu(t^{n+1/2})$ is sufficient to
maintain overall second order accuracy.  We compute
$\bar{\fatu} = \fatu^{n+1/2}$ using a first order two step
subscheme with time step $\Delta t/2$, as follows.

\medskip
\noindent\emph{Advection sub step} (with time step $\Delta t/2$,
advecting velocity $\fatu^n$): apply SSPRK3 to the advection
equations to obtain $(\fatu^*, D^*)$.

\medskip
\noindent\emph{Projection sub step}: solve
for $(\mathbf{u}^{n+1/2},D^{n+1/2})$ satisfying
\begin{gather}
    \mathbf{u}^{n+1/2} - \mathbf{u}^*
    + \frac{\Delta t}{2} f (\mathbf{u}^{n+1/2})^{\perp}
    = - g\, \frac{\Delta t}{2} \, \nabla (D^{n+1/2}+b),\\[4pt]
    D^{n+1/2} - D^*
    + \frac{\Delta t}{2}\nabla\cdot(\mathbf{u}^{n+1/2}\Bar{D})
    = 0.
\end{gather}
The result $\bar{\fatu} = \fatu^{n+1/2}$ is a first order
accurate approximation to $\fatu(t^{n+1/2})$.  To see that
using this advecting velocity in the second order scheme preserves
second order accuracy, consider the following system with frozen advecting velocity
\begin{align}
    \frac{\partial \mathbf{u}}{\partial t}
        + (\Bar{\fatu} \cdot\nabla)\fatu
        + f\mathbf{u}^{\perp}
        &= -g\nabla (D+b),\\
    \frac{\partial D}{\partial t}
        + \nabla \cdot (\Bar{\fatu} D)
        &= 0.
\end{align}
A Taylor expansion gives
$\fatu(t^{n+1}) = \fatu(t^n) + \Delta t\,\dot{\fatu}(t^{n+1/2})
+ \mathcal{O}(\Delta t^3)$.
Evaluating the right hand side at $t^{n+1/2}$ using the
first order approximation $\bar{\fatu}$ introduces only an
$\mathcal{O}(\Delta t^2)$ perturbation, so the overall
truncation error remains $\mathcal{O}(\Delta t^2)$.

In particular, $\bar{\fatu}\in V_h^1 \subset \Hdiv$, so that
the upwind flux in the DG advection step is well-defined and
the normal components of the advecting velocity are
single valued across element facets.

\subsection{Spatial discretisation}
\label{subsec:space-disc}

For a given triangulation $\Th$, let $V_h^1 \subset \Hdiv$ and $V_h^2 \subset L^2$ be
compatible finite element spaces, i.e. $\nabla\cdot$ maps $V_h^1$ onto $V_h^2$ with the image being the subspace of mean zero functions in
$V_h^2$.
Let $\hat{V}_h^1 \subset (L^2)^d$ denote the broken space
obtained by disassociating the degrees of freedom of $V_h^1$
on element facets; thus $V_h^1 \subset \hat{V}_h^1$.
Multiplying the momentum equation by test functions
$\fatw \in V_h^1$ (or $\fatw \in \hat{V}_h^1$ in the
advection step) and the continuity equation by
$\phi \in V_h^2$, integrating by parts element wise,
and summing over all elements gives the following fully
discrete equations.

\medskip
\noindent\textbf{Step 1 — Explicit step.}
Find $(\hat{\fatu}, \hat{D})\in V_h^1 \times V_h^2$ such that
\begin{align}
\big(\hat{\mathbf{u}} - \mathbf{u}^n,\, \fatw \big)_{\Th}
+ \frac{\Delta t}{2} \Big[
    \big( f(\mathbf{u}^n)^{\perp},\, \fatw \big)_{\Th}
  - g\,\big(D^n+b,\, \nabla_h \cdot \fatw\big)_{\Th}
\Big] &= 0,
\label{eq:disc_exp_u}\\[6pt]
\big(\hat{D} - D^n,\, \phi \big)_{\Th}
+ \frac{\Delta t}{2}
  \big(\nabla_h\cdot(\mathbf{u}^n\Bar{D}),\, \phi \big)_{\Th}
&= 0,
\label{eq:disc_exp_D}
\end{align}
for all $\fatw \in V_h^1$ and $\phi \in V_h^2$, respectively.
Here $\nabla_h\cdot$ and $\nabla_h$ denote the element wise
(broken) divergence and gradient, respectively.

\medskip
\noindent\textbf{Step 2 — DG-upwind advection step.}
We define the DG bilinear form for velocity advection.
For a given advecting velocity $\Bar{\fatu}\in \Hdiv$, define
\begin{equation}
\label{eq:L_u}
    \mathcal{L}_{\Bar{\mathbf{u}}}^{\fatu}(\fatw, \fatu)
    :=  \int_{\Th} \nabla_h \cdot (\fatw \otimes \Bar{\fatu})
        \cdot \fatu\, \dx
      - \int_{\Eh} (\Bar{\fatu}_+ \cdot \fatn_+) \,
        \widetilde{\fatu} \cdot\llbracket \fatw \rrbracket_+
        \, \dS + C_{\mathrm{sph}},
\end{equation}
where $\Eh$ is the set of facets in $\Th$, $\llbracket \fatw \rrbracket_+:= \fatw_+ - \fatw_-$
denotes the jump across a cell facet, $\fatn^{\pm}$ are the
outward normals from the $\pm$ sides, and
\begin{equation}
    C_{\mathrm{sph}}
    = \int_{\Eh}\left( \Bar{\fatu}^+ \cdot \fatn^+ \right)
      \left( \widetilde{\fatu} \cdot \widetilde{\fatn} \right)
      \left( \underline{\fatw}
             \left[ \fatn^+ + \fatn^- \right]\right) \dS.
\end{equation}
Since $\Bar{\fatu}\in \Hdiv$ the normal component is
single valued: $\Bar{\fatu}_+\cdot\fatn_+ =
-\Bar{\fatu}_-\cdot\fatn_-$.  Following \cite{bernard2009high}, the term $C_{\mathrm{sph}}$ is a
spherical-geometry correction ensuring that both sides of the
jump $\llbracket \fatw \rrbracket_+$ are evaluated in the
tangent space of the upwind value, where the upwind velocity is
\begin{equation}
    \widetilde{\fatu} = \begin{cases}
        \fatu^+ & \text{if }\fatu^+ \cdot \fatn^+ \geq 0,\\
        \fatu^- & \text{if }\fatu^+ \cdot \fatn^+ < 0,
    \end{cases}
\end{equation}
and $\underline{\fatw}$ denotes the downwind test function value.
For depth advection, the standard upwind DG form is
\begin{equation}
\label{eq:L_D}
    \mathcal{L}_{\Bar{\mathbf{u}}}^D(\phi, D)
    :=  \int_{\Th} D \, \nabla_h \phi \cdot \Bar{\mathbf{u}}
        \, \dx
      - \int_{\Eh} \widetilde{D} \,
        \llbracket \phi \, \Bar{\mathbf{u}} \rrbracket \dS,
\end{equation}
where $\widetilde{D}$ is the upwind depth value and
$\llbracket\cdot\rrbracket$ the normal jump.

The SSPRK3 advection step reads: find
$(\fatu^*, D^*)\in \hat{V}_h^1 \times V_h^2$ such that
\begin{align}
  \big(\mathbf{u}^{(1)},\, \mathbf{w}\big)_{\Th}
  &= \big( \hat{\fatu},\, \mathbf{w}\big)_{\Th}
   + \Delta t \,\mathcal{L}_{\Bar{\mathbf{u}}}^{\fatu}
     (\mathbf{w}, \hat{\fatu}),
  \label{eq:ssprk1_u}\\
  \big(\fatu^{(2)},\, \mathbf{w}\big)_{\Th}
  &= \tfrac{3}{4}\big( \hat{\fatu},\, \mathbf{w}\big)_{\Th}
   + \tfrac{1}{4}\Big[
     \big(\fatu^{(1)},\, \fatw \big)_\Th
   + \Delta t\,\mathcal{L}_{\Bar{\mathbf{u}}}^{\fatu}
     (\mathbf{w}, \fatu^{(1)}) \Big],
  \label{eq:ssprk2_u}\\
  \big(\fatu^{*},\, \mathbf{w}\big)_{\Th}
  &= \tfrac{1}{3}\big( \hat{\fatu},\, \mathbf{w}\big)_{\Th}
   + \tfrac{2}{3}\Big[
     \big(\fatu^{(2)},\, \fatw \big)_\Th
   + \Delta t\,\mathcal{L}_{\Bar{\mathbf{u}}}^{\fatu}
     (\mathbf{w}, \fatu^{(2)}) \Big],
  \label{eq:ssprk3_u}
\end{align}
for all $\mathbf{w} \in \hat{V}_h^1$, and

\begin{align}
  \big(D^{(1)},\, \phi\big)_{\Th}
  &= \big(\hat{D},\, \phi\big)_{\Th}
   + \Delta t\,\mathcal{L}_{\Bar{\mathbf{u}}}^D
     (\phi, \hat{D}-\Bar{D}),
  \label{eq:ssprk1_D}\\
  \big(D^{(2)},\, \phi\big)_{\Th}
  &= \tfrac{3}{4}\big(\hat{D},\, \phi\big)_{\Th}
   + \tfrac{1}{4}\Big[
     \big(D^{(1)},\, \phi\big)_\Th
   + \Delta t\,\mathcal{L}_{\Bar{\mathbf{u}}}^D
     (\phi, D^{(1)}-\Bar{D}) \Big],
  \label{eq:ssprk2_D}\\
  \big(D^{*},\, \phi\big)_{\Th}
  &= \tfrac{1}{3}\big(\hat{D},\, \phi\big)_{\Th}
   + \tfrac{2}{3}\Big[
     \big(D^{(2)},\, \phi\big)_\Th
   + \Delta t\,\mathcal{L}_{\Bar{\mathbf{u}}}^D
     (\phi, D^{(2)}-\Bar{D}) \Big],
  \label{eq:ssprk3_D}
\end{align}
for all $\phi \in V_h^2$.

\medskip
\noindent\textbf{Step 3 — Projection step.}
Find $(\fatu^{n+1}, D^{n+1})\in V_h^1 \times V_h^2$ such that
\begin{align}
\big(\mathbf{u}^{n+1} - \mathbf{u}^*,\, \mathbf{w}\big)_{\Th}
+ \frac{\Delta t}{2} \Big[
    \big(f(\mathbf{u}^{n+1})^{\perp},\, \mathbf{w}\big)_\Th
  - g\big(D^{n+1}+b,\, \nabla_h \cdot \mathbf{w} \big)_\Th
\Big] &= 0,
\label{eq:disc_proj_u}\\[6pt]
\big(D^{n+1} - D^*,\, \phi\big)_{\Th}
+ \frac{\Delta t}{2}
  \big(\nabla_h \cdot (\mathbf{u}^{n+1}\Bar{D}),\, \phi \big)_\Th
&= 0,
\label{eq:disc_proj_D}
\end{align}
for all $\fatw \in V_h^1$ and $\phi \in V_h^2$, respectively.

\subsection{Finite element spaces}
\label{subsec:fe-spaces}

All results in Section \ref{sec:results} use the second order
compatible pair
\begin{equation}
    V_h^1 = BDM_2,\qquad V_h^2 = DG_1,
\end{equation}
where $BDM_2$ is the second order Brezzi-Douglas-Marini
element \cite{brezzi1985two} and $DG_1$ denotes the space of
discontinuous piecewise linear polynomials.  These spaces satisfy
$\nabla\cdot BDM_2 = DG_1$, preserving the discrete de Rham
structure. The spaces are defined on an icosahedral triangulation of the sphere, and the $\Hdiv$ basis functions are mapped to each cell via the contravariant Piola transform, so that discrete velocities are tangent to the mesh surface. In the advection step the broken space
$\hat{V}_h^1 = BDM_2^d$ is used, obtained by breaking $BDM_2$
on element interfaces.  With this pair, second order spatial
convergence $\mathcal{O}(h^2)$ is expected and confirmed in
Section \ref{subsec:convergence}.

\subsection{Solver strategy}
\label{subsec:solvers}

Scalable, parameter robust solvers are identified for each step.
Throughout, iteration counts are verified to be independent of
$\Delta t$ and $\Delta x$.  The spaces and solvers are
summarised in Table \ref{tab:solvers}.

\begin{table}[H]
    \centering
    \small
    \begin{tabular}{lcccp{5cm}}
        \hline
        Step & $\fatu$-space & $D$-space & Coupled? & Solver \\
        \hline
        Explicit   & $BDM_2$   & $DG_1$ & No  &
            GMRES + block ILU ($\fatu$);
            cellwise direct LU ($D$)\\[4pt]
        Advection  & $BDM_2^d$ & $DG_1$ & No  &
            cellwise direct LU \\[4pt]
        Projection & $BDM_2$   & $DG_1$ & Yes &
            hybridisation + direct LU \\
        \hline
    \end{tabular}
    \caption{Finite element spaces and parameter robust solvers
             for each step of the time stepping scheme.}
    \label{tab:solvers}
\end{table}

\paragraph{Explicit step}
The velocity and depth equations \eqref{eq:disc_exp_u}--\eqref{eq:disc_exp_D}
are decoupled.  The depth mass matrix is cellwise block diagonal and can be efficiently inverted by 
 direct solvers in parallel. The solve factorisations are cached and reused for each timestep.
The velocity mass matrix is globally coupled but can be
solved by GMRES with block ILU preconditioning.  The GMRES
iteration count is constant and independent of $\Delta t$ and
$\Delta x$. For example, with a relative tolerance of $10^{-7}$, the count
for the Galewsky et al test is uniformly equal to 7 across all mesh
levels and time steps tested.

\paragraph{Advection step}
The left hand side of each SSPRK3
stage \eqref{eq:ssprk1_u}--\eqref{eq:ssprk3_D} is a mass
matrix.  Since $\fatu\in\hat{V}_h^1$ and $D\in V_h^2$ are
both in discontinuous spaces, the mass matrices are
block diagonal.  Each stage is solved by cellwise direct solvers in parallel, cached and reused.

\paragraph{Projection step}
The coupled wave
system \eqref{eq:disc_proj_u}--\eqref{eq:disc_proj_D} is
solved via hybridisation \cite{cockburn2009unified}.
Hybridisation introduces trace variables (Lagrange multipliers)
supported on element facets and breaks the global $\Hdiv$
velocity space, transforming the coupled system into a smaller,
sparse condensed system for the traces alone.  The remaining
degrees of freedom are recovered by local elementwise
postprocessing.  For the present two dimensional problem the
condensed system is solved by a direct LU factorisation. The LU
factors are cached and reused, since the equation does not 
change between timesteps.
In
a three dimensional model, this global solve could be replaced
by a multigrid method, preserving the scalable structure of
the algorithm.

%% ============================================================
%%  Section 4 — Numerical Results
%% ============================================================
\section{Numerical Results}
\label{sec:results}

We assess the proposed scheme through three standard benchmarks
for global shallow water models: the Williamson 5 test
(topographic forcing and Rossby wave generation), the Galewsky et al
barotropic instability test (balanced flow and slow-fast mode
separation), and the L\"{a}uter unsteady solid body rotation
test (quantitative convergence analysis against an exact
solution).  All computations use the finite element pair
$(BDM_2,\,DG_1)$ described in Section \ref{subsec:fe-spaces},
implemented in the Firedrake finite element
framework \cite{FiredrakeUserManual} on icosahedral sphere
meshes.

As {noted} in Section~\ref{subsec:subcycling}, the wave
part forms an unconditionally stable Crank--Nicolson scheme,
so stability is governed solely by the CFL condition of
the upwind DG advection step.  All quoted Courant numbers are the
Courant number $C$~\eqref{eq:courant} of the full time step
$\Delta t$. In all experiments reported here the scheme was stable
with $C/N = 0.5$ in the advection substeps, i.e.\ $C = 0.5$ without
subcycling and $C = 1$ with $N = 2$ subcycles.

\subsection{Williamson 5: Zonal flow over an isolated mountain}
\label{subsec:W5}

The Williamson 5 test \cite{williamson1992standard} is a widely
used benchmark for topographic forcing and conservation
properties.  It initialises a geostrophically balanced zonal
flow that is subsequently perturbed by an isolated conical
mountain.  The test exercises the scheme's ability to generate
and sustain large scale Rossby waves while conserving energy and
enstrophy over extended integration periods.

\subsubsection{Initial conditions and parameters}
The initial velocity (in Cartesian coordinates) and
geostrophically balanced total height are
\begin{equation}
    \fatu = \frac{u_0}{R_0}
    \begin{pmatrix}-y \\ x \\ 0\end{pmatrix},\qquad
    D + b = H - \frac{1}{g}
      \!\left(R_0\Omega u_0 + \frac{u_0^2}{2}\right)
      \frac{z^2}{R_0^2},
\end{equation}
with $u_0 = 20$ m\,s$^{-1}$, $R_0 = 6.371\times10^6$ m,
$\Omega = 7.292\times10^{-5}$ s$^{-1}$, $H = 5960$ m.
The conical mountain topography (in spherical coordinates) is
\begin{equation}
    b = b_0\!\left(1 - r/R_m\right),\qquad
    r^2 = \min\!\left\{R_m^2,\;
           (\lambda-\lambda_c)^2 + (\theta-\theta_c)^2\right\},
\end{equation}
with $b_0 = 2000$ m, $R_m = \pi/9$, $\lambda_c = -\pi/2$,
$\theta_c = \pi/6$.

\subsubsection{Results and stability diagnostics}
\begin{figure}[htbp]
    \centering
    \includegraphics[width=\textwidth]%
        {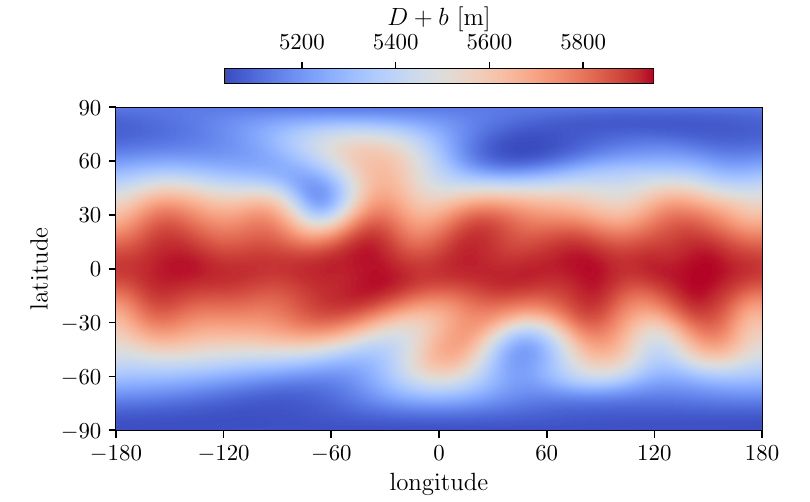}

    \vspace{0.5em}
    \includegraphics[width=\textwidth]%
        {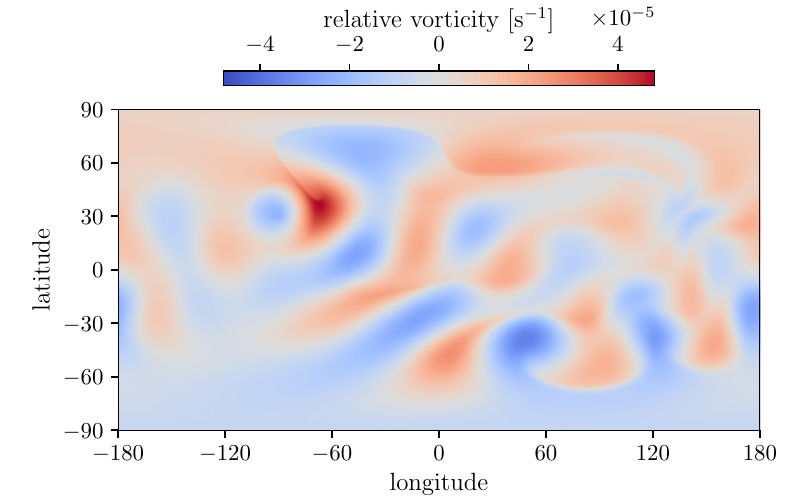}
    \caption{Williamson 5 test at $t=15$ d: total height $D + b$ (top)
    and relative vorticity (bottom). Icosahedral mesh level 6,
    $\Delta t = 1200$ s, $N=2$ advection subcycles ($C=1$).}
    \label{fig:W5}
\end{figure}

The test was run on an icosahedral mesh of refinement level 6
with $\Delta t = 1200$ s and $N = 2$
advection subcycles, giving $C = 1$ at
$t = 15$ days.  Figure \ref{fig:W5} shows the total height $D + b$
and the relative vorticity at $t = 15$ days.  The expected slowly
propagating Rossby waves are visible downstream of the
mountain, and the vorticity field is smooth with no spurious
grid scale oscillations.

To quantify the stability properties of the scheme, we
monitor three diagnostics.  The relative energy error is
\begin{equation}
    e_E(t) = \frac{E_t - E_0}{E_0},\quad
    E_t = \frac{1}{2}\!\int_\Omega\!\Big(
          D(t)|\fatu(t)|^2
          + g\,(D(t)+b)^2\Big)\dx,
\end{equation}
where $E_t$ omits the time-independent term
$-\frac{g}{2}\int_\Omega b^2\dx$ of the total energy, which
does not affect $E_t - E_0$;
the relative potential enstrophy error is
\begin{equation}
    e_Z(t) = \frac{\hat{Z}_t - \hat{Z}_0}{\hat{Z}_0},\quad
    \hat{Z}_t = \frac{1}{2}\!\int_\Omega q(t)^2 D\,\drm x,
\end{equation}
where $q$ is the potential vorticity, and the change in the root-mean-square divergence is
\begin{equation}
    R_{\mathrm{div}}(t) =
    \frac{\|\nabla\cdot\mathbf{u}(t)\|_{L^2}
          - \|\nabla\cdot\mathbf{u}(0)\|_{L^2}}
         {|\Omega|^{1/2}},
    \qquad |\Omega| = \int_\Omega 1\,\dx.
\end{equation}
%We normalise by $|\Omega|^{1/2}$ rather than by
%$\|\nabla\cdot\mathbf{u}(0)\|_{L^2}$, since the initial flow is
%nearly non-divergent, so the latter is close to zero.
$R_{\mathrm{div}}$ serves as an
indicator of the effectiveness of the projection step in
controlling spurious divergence modes.

\begin{figure}[htbp]
    \centering
    \includegraphics[width=0.9\textwidth]%
        {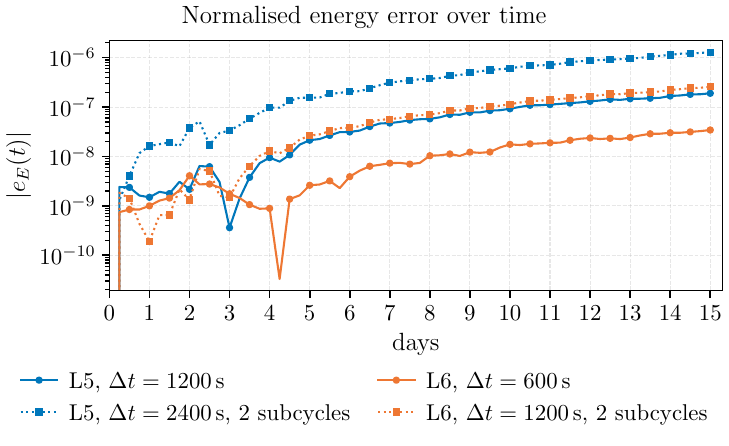}
        \caption{Williamson 5: relative energy error
    $|e_E|$ (logarithmic scale) over 15 days.  For the larger time
    step ($C=1$), the advection step is subcycled twice ($N=2$), so that
    each advection substep has the same Courant number as the smaller time step
    ($C=0.5$).  Mesh level 5: $\Delta t=1200$ s and $2400$ s;
    mesh level 6: $\Delta t=600$ s and $1200$ s.}

    \label{fig:W5_energy_error}
\end{figure}

\begin{figure}[htbp]
    \centering
    \includegraphics[width=0.9\textwidth]%
        {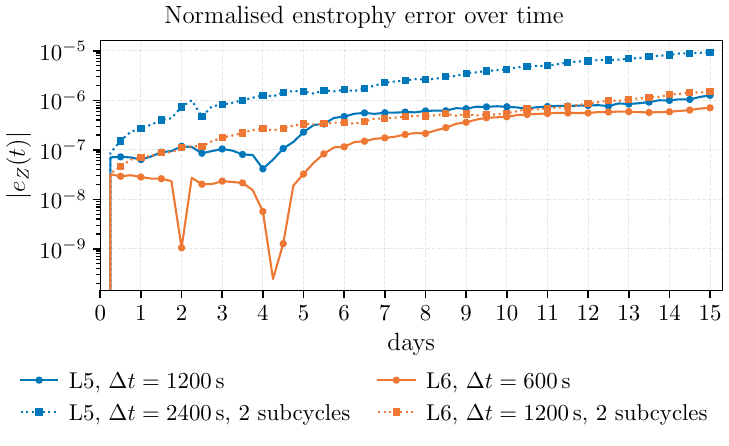}
    \caption{Williamson 5: relative potential enstrophy error
    (absolute value, logarithmic scale)
    over 15 days.  Run parameters as in
    Figure \ref{fig:W5_energy_error}.}
    \label{fig:W5_enstrophy_error}
\end{figure}

\begin{figure}[htbp]
    \centering
    \includegraphics[width=0.9\textwidth]%
        {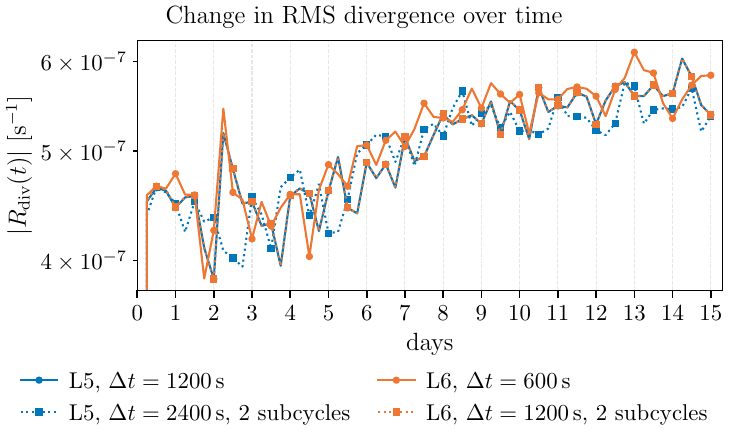}
    \caption{Williamson 5: change in the RMS
    divergence $|R_{\mathrm{div}}|$ (logarithmic scale)
    over 15 days.  Run parameters as
    in Figure \ref{fig:W5_energy_error}.}
    \label{fig:W5_divl2_error}
\end{figure}

Figures \ref{fig:W5_energy_error}--\ref{fig:W5_divl2_error}
show these diagnostics for mesh levels 5 and 6, with and
without subcycling, over 15 days.  Energy errors are slightly
increasing but remain bounded throughout at order
$10^{-8}$--$10^{-6}$; subcycling leads to errors approximately
one order of magnitude larger, but still bounded.  Enstrophy
errors are similarly bounded at order $10^{-5}$--$10^{-6}$.
$R_{\mathrm{div}}$ is bounded throughout at order
$10^{-7}$\,s$^{-1}$, confirming that the $\Hdiv$-projection step
effectively controls divergent modes.

\subsection{Galewsky et al barotropic instability test}
\label{subsec:Galewsky}

The Galewsky et al test \cite{galewsky2004initial} provides a test of
the scheme's ability to maintain the
separation between fast gravity waves and balanced geostrophic
motion.  It initialises a barotropically unstable midlatitude
jet and adds a small perturbation that triggers
instability.  The subsequent evolution passes through
two qualitatively distinct phases: a fast initial adjustment
dominated by gravity wave radiation, followed by a slow
roll up of the vorticity field that is characteristic of
barotropic instability.

\subsubsection{Initial conditions and parameters}
The initial zonal velocity (in spherical coordinates) is
\begin{equation}
    \fatu(\phi) = \begin{cases}
        0
        & \phi \leq \phi_0,\\[4pt]
        \dfrac{u_{\max}}{e_n}
        \exp\!\left[\dfrac{1}{(\phi-\phi_0)(\phi-\phi_1)}
             \right]
        & \phi_0 < \phi < \phi_1,\\[4pt]
        0
        & \phi \geq \phi_1,
    \end{cases}
\end{equation}
with $u_{\max} = 80$ m\,s$^{-1}$, $\phi_0 = \pi/7$,
$\phi_1 = \pi/2 - \phi_0$, and
$e_n = \exp(-4/(\phi_0-\phi_1)^2)$.  A balanced depth field
$D_0(\phi)$ is obtained by integrating the balance equation and shifted so that its global mean is $H_0 = 10$ km.  
A perturbation
\begin{equation}
    h'(\lambda,\phi) = \hat{h}\cos(\phi)\,
    \mathrm{e}^{-(\lambda/\alpha)^2}\,
    \mathrm{e}^{-[(\phi_2-\phi)/\beta]^2},
    \quad -\pi\leq\lambda\leq\pi,
\end{equation}
centred at $\phi_2 = \pi/4$ is added to $D_0$ to initiate the
instability.  All parameters follow \cite{galewsky2004initial}.
The reference depth used in the splitting is
$\bar{D} = H_0 = 10$ km, equal to the global mean of $D_0$. 

\subsubsection{Results and stability diagnostics}
\begin{figure}[htbp]
    \centering
    \includegraphics[width=\textwidth]%
        {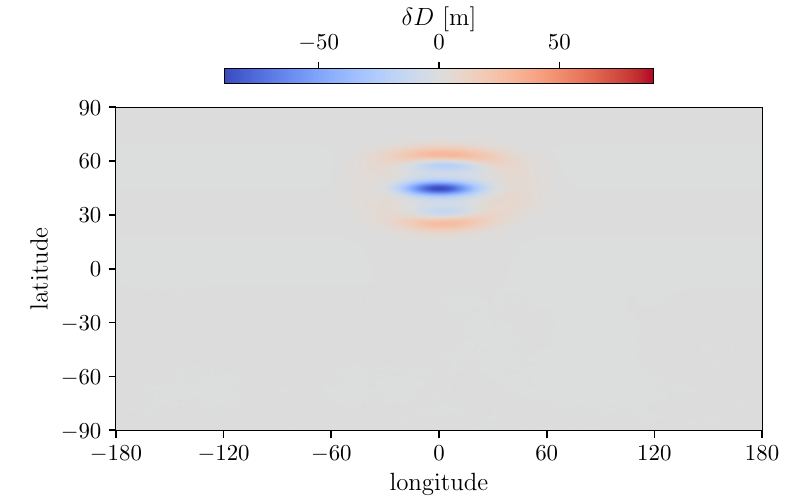}

    \vspace{0.5em}
    \includegraphics[width=\textwidth]%
        {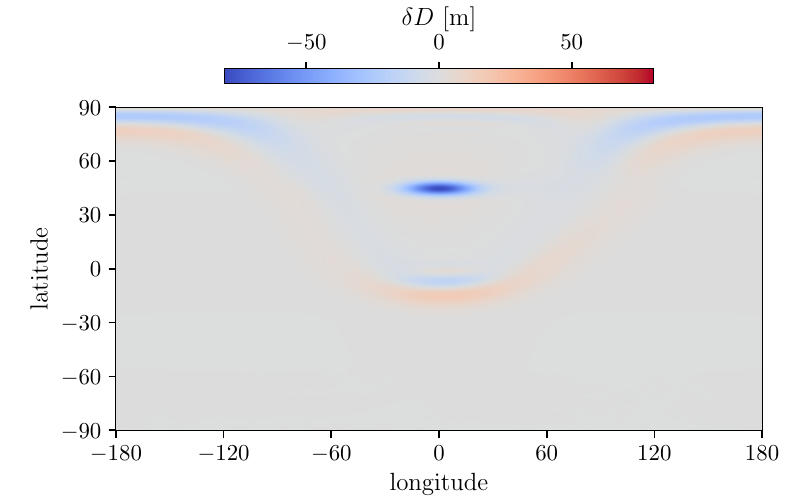}
    \caption{Galewsky et al test: fast-wave adjustment phase, depth
    perturbation $\delta D = D(t)-D(0)$ at the start ($t=2$ h, top)
    and end ($t=6$ h, bottom) of the adjustment. Mesh level 6,
    $\Delta t = 550$ s, $N=2$ subcycles ($C=1$).}
    \label{fig:Galewsky_D}
\end{figure}

In this test, we used an icosahedral mesh of refinement level 6
with $\Delta t = 550$ s and $N = 2$ subcycles ($C = 1$ at $t=6$ d).
Figure \ref{fig:Galewsky_D} shows the depth perturbation $\delta D = D(t) - D(0)$, i.e.\ the deviation
from the initial perturbed depth $D_0 + h'$, at $t = 2$ and $6$ hours, illustrating the fast gravity wave adjustment phase.

\begin{figure}[htbp]
    \centering
    \includegraphics[width=\textwidth]%
        {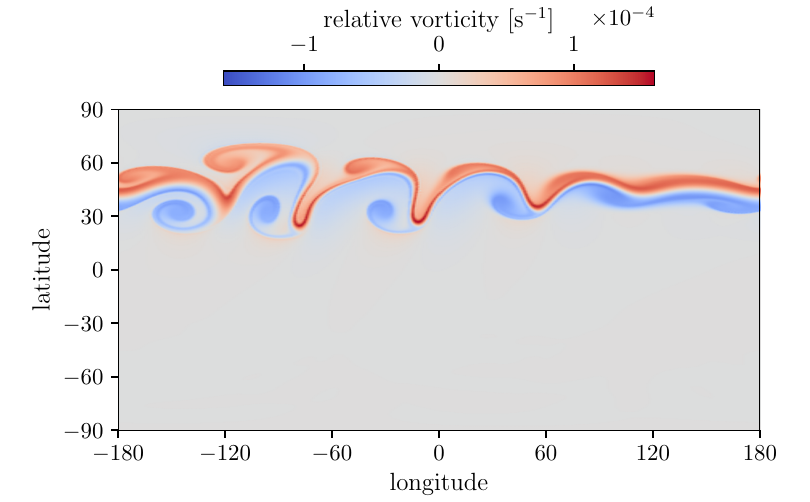}
    \caption{Galewsky et al test: relative vorticity at $t=144$ h
    ($t=6$ d), showing the fully developed barotropic
    instability. Same run parameters as
    Figure \ref{fig:Galewsky_D}.}
    \label{fig:Galewsky_vort}
\end{figure}

Figure \ref{fig:Galewsky_vort} shows the relative vorticity
at $t = 144$ hours, exhibiting the characteristic vorticity
roll-up consistent with the reference results
of \cite{galewsky2004initial}.  The clean separation of the
fast and slow phases confirms that the scheme accurately
represents both balanced and unbalanced dynamics without
spurious mode coupling.

\begin{figure}[htbp]
    \centering
    \includegraphics[width=0.9\textwidth]%
        {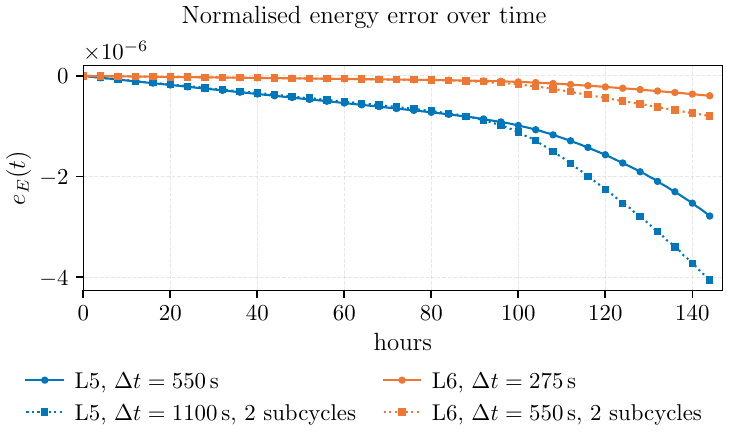}
        \caption{Galewsky et al test: relative energy error over 144 h.
    For the larger time step ($C=1$), the advection step is subcycled
    twice ($N=2$); the smaller time step has $C=0.5$.  Mesh level 5:
    $\Delta t=550$ s and $1100$ s; mesh level 6: $\Delta t=275$ s
    and $550$ s.}

    \label{fig:Galewsky_energy_error}
\end{figure}

\begin{figure}[htbp]
    \centering
    \includegraphics[width=0.9\textwidth]%
        {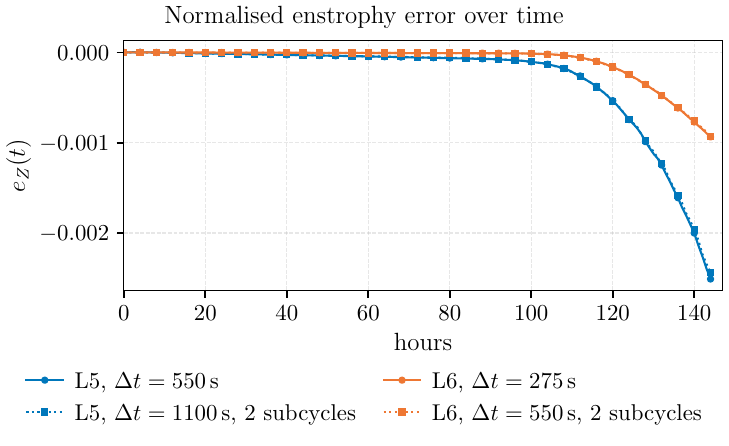}
    \caption{Galewsky et al test: relative potential enstrophy error
    over 144 h.  Run parameters as in
    Figure \ref{fig:Galewsky_energy_error}.}
    \label{fig:Galewsky_enstrophy_error}
\end{figure}

\begin{figure}[htbp]
    \centering
    \includegraphics[width=0.9\textwidth]%
        {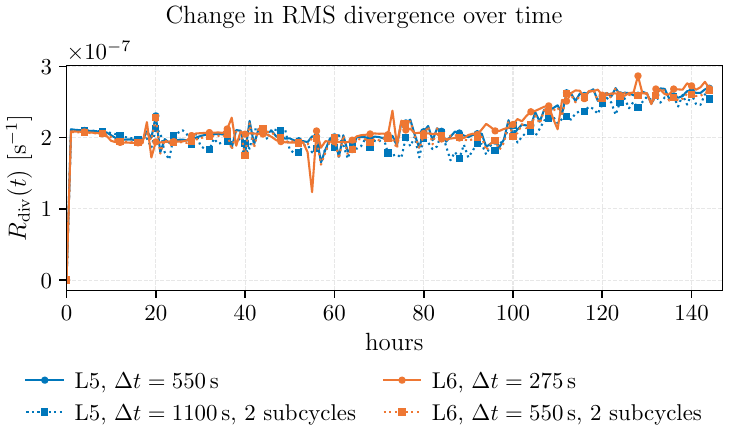}
    \caption{Galewsky et al test: change in the RMS
    divergence $R_{\mathrm{div}}$ over 144 h. Run parameters as in
    Figure \ref{fig:Galewsky_energy_error}.}
    \label{fig:Galewsky_divl2_error}
\end{figure}

Figures \ref{fig:Galewsky_energy_error}--\ref{fig:Galewsky_divl2_error}
show energy, enstrophy, and divergence diagnostics for mesh
levels 5 ($\Delta t = 550$ s and $1100$ s) and 6 ($\Delta t = 275$ s and $550$ s).
In contrast to the Williamson 5 results, the Galewsky et al test
exhibits a small energy damping rather than a slow energy drift.
This is consistent with the inherent numerical dissipation of
the upwind DG scheme and is more pronounced here because the
Galewsky et al test is more strongly advection dominated.
Importantly, the damping rate decreases with mesh refinement:
the level 6 results show noticeably less dissipation than
level 5, consistent with numerical convergence.
$R_{\mathrm{div}}$ remains bounded at order
$10^{-7}$\,s$^{-1}$ throughout.

\subsection{Convergence analysis: unsteady solid-body rotation}
\label{subsec:convergence}

L\"{a}uter et al. \cite{lauter2005unsteady} derived initial
conditions for the spherical shallow water equations for which
exact analytical solutions are available.  The solution is a
solid-body rotation whose axis is tilted with respect to the
rotation axis of the sphere, causing the fields to rotate
around the Earth and providing a genuinely unsteady, non-trivial
test for convergence analysis.

\subsubsection{Exact solution}
Let $\mathbf{b}_i(t)$ be the coordinate system rotating about
the $z$-axis with angular velocity $\Omega$:
\begin{equation}
    \mathbf{b}_1(t) = \cos(\Omega t)\,\mathbf{e}_1
                    + \sin(\Omega t)\,\mathbf{e}_2,\quad
    \mathbf{b}_2(t) = -\sin(\Omega t)\,\mathbf{e}_1
                    + \cos(\Omega t)\,\mathbf{e}_2,\quad
    \mathbf{b}_3(t) = \mathbf{e}_3,
\end{equation}
and define the coordinate transform
$\varphi_t(\mathbf{x}) = (\mathbf{x}\cdot\mathbf{b}_i)\mathbf{e}_i$.
With rotation axis $\mathbf{c} = -\sin(\alpha)\,\mathbf{e}_1
+ \cos(\alpha)\,\mathbf{e}_3$ and
$\alpha = \pi/4$, the exact velocity and total height $D + b$ are
\begin{align}
    \mathbf{u}(\mathbf{x}, t)
        &= u_0\,\varphi_t(\mathbf{c}) \times \mathbf{n},\\
    D(\mathbf{x}, t) + b(\mathbf{x})
        &= \frac{1}{g}\left(k_1
           - \frac{(u_0\,\varphi_t(\mathbf{c})\cdot\mathbf{n}
             + \boldsymbol{\Omega}\cdot\mathbf{x})^2}{2}
           + \frac{(\boldsymbol{\Omega}\cdot\mathbf{x})^2}{2}\right),
\end{align}
with bottom topography $b(\mathbf{x}) = \frac{1}{2g}(\boldsymbol{\Omega}\cdot\mathbf{x})^2$.  We use $u_0 = 2\pi R_0/(12\ \mathrm{d})$ and
$k_1 = 133681$ m$^2$\,s$^{-2}$.

\subsubsection{Convergence results}
\begin{figure}[htbp]
    \centering
    \includegraphics[width=0.95\linewidth]%
        {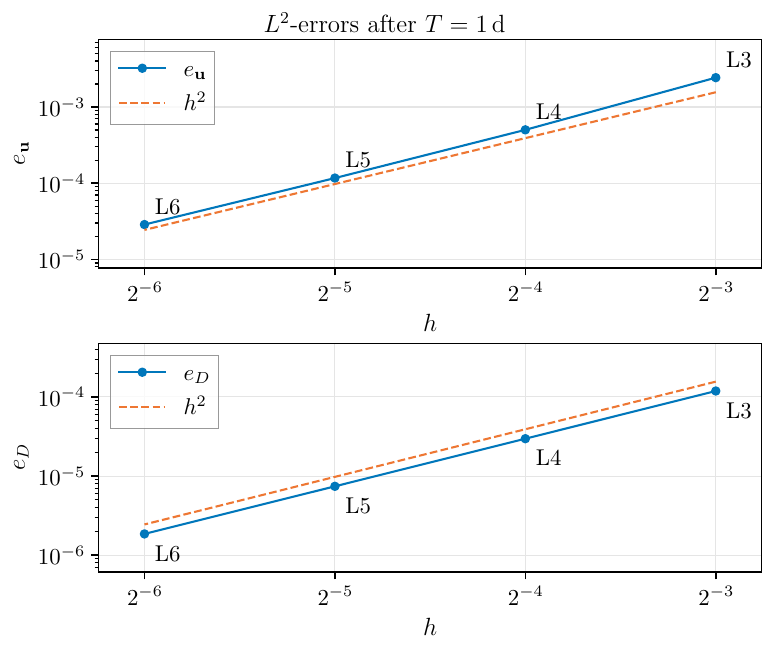}
    \caption{Log-log convergence plot: normalised $L^2$-errors
    for $\fatu$ (top) and $D+b$ (bottom) at $T=1$ day, run with
    constant $\Delta t/\Delta x$ on mesh refinement levels 3--6
    (L3--L6).  The reference slope $h^2$ is shown for
    comparison.  Quadratic convergence is observed for both
    fields.}
    \label{fig:convergence}
\end{figure}

The normalised $L^2$-errors are defined as
\begin{equation}
    e_{\fatu} = \frac{\|\fatu_h - \fatu_{\rm exact}\|_{L^2}}
                     {\|\fatu_{\rm exact}\|_{L^2}},
    \qquad
    e_D = \frac{\|(D_h + b_h) - (D+b)_{\rm exact}\|_{L^2}}
               {\|(D+b)_{\rm exact}\|_{L^2}}.
\end{equation}
The scheme was run to $T = 1$ day on icosahedral meshes of
refinement levels 3--6 with time-steps
$\Delta t = 400, 200, 100, 50$ s respectively, maintaining a
constant $\Delta t/\Delta x$ ratio so that spatial and temporal
errors converge together.  Figure \ref{fig:convergence} shows
log-log plots of $e_{\fatu}$ and $e_D$, confirming quadratic
convergence $\mathcal{O}(h^2)$ for both fields in agreement
with the expected second-order accuracy of the $(BDM_2, DG_1)$
spatial pair and the symmetric time-stepping scheme.

\subsection{Parallel performance}
\label{subsec:timings}

To briefly touch upon the practicality of the scheme,
table \ref{tab:timings} reports wall-clock times for the
Williamson 5 test ($T = 15$ d, mesh level 6,
$\Delta t = 1200$ s, $N=2$ subcycles) run in parallel on a
dual Intel Xeon Silver 4215 system (8 cores per processor,
2.50 GHz) across 6, 12, and 15 cores.  The projection step
dominates cost ($\approx 60\%$ of total), followed by the
advection step ($\approx 30\%$) and the explicit step
($< 10\%$).  The advection step exhibits super-linear
speed-up between 6 and 12 cores, likely due to improved
cache usage and domain-decomposition effects with the
default Firedrake mesh partitioner; this could be further
improved by tuning partitioner options.  The projection
step scales well, reflecting the reduction of the global
system to a smaller condensed problem via hybridisation.
Overall, the scheme achieves a speed-up factor of
approximately 2 going from 6 to 12 cores, with further
modest improvement at 15 cores. This is fairly typical
single node performance.

\begin{table}[H]
    \centering
    \begin{tabular}{lrrr}
        \hline
        \textbf{Step} & \textbf{6 cores}
                      & \textbf{12 cores}
                      & \textbf{15 cores} \\
        \hline
        Total         & 5.82 & 3.25 & 2.88 \\
        Explicit step & 0.29\ \ (5.0\%) & 0.24\ \ (7.2\%) & 0.25\ \ (8.7\%) \\
        Advection step& 1.85\ (31.8\%) & 0.83\ (25.5\%) & 0.78\ (27.2\%) \\
        Projection step& 3.57\ (61.4\%) & 2.10\ (64.6\%) & 1.77\ (61.4\%) \\
        \hline
    \end{tabular}
    \caption{Wall-clock times (in units of $10^3$ s) for the
    Williamson 5 test across 6, 12, and 15 CPU cores.
    Percentages indicate the fraction of total run time for
    each step.}
    \label{tab:timings}
\end{table}

%% ============================================================
%%  Section 5 — Conclusion
%% ============================================================
\section{Conclusion}
\label{sec:conclusion}

We have presented and analysed a second order, semi-implicit
time stepping scheme for the shallow water equations on the
sphere, discretised in space with compatible finite elements.
The scheme is a symmetric three step splitting: a forward Euler
half step for the wave part using $\Hdiv$-conforming velocities,
a full SSPRK3 upwind DG advection step using a broken $\Hdiv$
space, and a backward Euler half step for the wave part that
simultaneously restores $\Hdiv$ conformity of the velocity via
an implicit projection.  The symmetry of this structure
guarantees global second order accuracy in time, and the wave
part (Steps 1 and 3 combined) forms an unconditionally stable
Crank--Nicolson scheme, confining the stability restriction
to the CFL condition of the upwind DG advection step alone.

Three standard benchmark tests were used to validate the
scheme.  The Williamson 5 test demonstrated stable, physically
correct development of large-scale Rossby waves in the presence
of topography over integration periods of 15 days, with
bounded energy errors of order $10^{-8}$--$10^{-6}$, enstrophy
errors of order $10^{-6}$--$10^{-5}$, and a change in the
root-mean-square divergence bounded at order $10^{-7}$\,s$^{-1}$,
confirming that the projection step effectively controls spurious
divergent modes.  The Galewsky et al barotropic instability test
confirmed that the scheme accurately reproduces both the fast
initial gravity-wave adjustment and the subsequent slow
vorticity roll-up without spurious mode coupling.  A small
energy damping was observed in the Galewsky et al test, consistent
with the inherent dissipation of the upwind DG scheme and
decreasing with mesh refinement, confirming its convergent
nature.  The L\"{a}uter unsteady solid-body rotation test
provided quantitative confirmation of second order spatial and
temporal convergence for both the velocity and the total height.
Across all tests, stable results were obtained for Courant
numbers up to $C = 1$ using second order finite elements, with
two advection subcycles, i.e.\ a Courant number of $0.5$ in each
advection substep.

An important practical feature of the scheme is the subcycling
strategy: by performing $N$ advection sub steps of size
$\Delta t/N$ within a single global time step $\Delta t$, the
Courant number of each advection substep is reduced to $C/N$,
allowing larger global time steps than the
bare CFL condition would permit, without incurring additional
implicit solves.  This is particularly advantageous because the
projection step — which involves a global implicit system
solved via hybridisation — is the most expensive component
($\approx 60\%$ of total run time), while the advection step
($\approx 30\%$) involves only block diagonal mass matrix
applications.  The hybridisation strategy also makes the solver
intrinsically scalable. By reducing the global system to a
smaller condensed problem on element facets, it is well suited
to both direct solves (as used here in 2D) and to multigrid
methods (appropriate for 3D).  Parallel timing results on up
to 15 cores confirm good scaling properties.

Compared to the finite volume semi-implicit method of
Benacchio and Klein \cite{benacchio2019semi}, the present scheme
operates natively within the compatible finite element framework,
and avoids nested solver loops.  The stable Courant number of
$C = 0.5$ achieved here without subcycling is broadly comparable
to the value of 1 reported in \cite{benacchio2019semi} for first
order finite volumes: the difference reflects the larger stencil of
second order DG spaces, and the effective degrees of freedom
per stable time step are similar.

Several directions for future work are natural.  The most
immediate is the extension of the scheme to the full
three-dimensional compressible Euler equations, which share the
same slow-fast split structure and for which the splitting
and solver components validated here are directly applicable.
Higher order finite element spaces ($p > 2$) are also of
interest: as $p$ increases, the ratio of facet to interior
degrees of freedom decreases, so the condensed hybridised
system grows more slowly than the overall problem size,
making higher order methods computationally attractive.
Finally, the scheme's stability threshold and energy properties
could be further improved by replacing SSPRK3 with
higher order invariant-domain-preserving
Runge-Kutta methods \cite{ern2022invariant}, or by
exploring alternative limiters for the DG advection step.

%% ============================================================
%%  Declaration of competing interests
%% ============================================================
\section*{Declaration of competing interests}
The authors declare that they have no known competing financial
interests or personal relationships that could have appeared to
influence the work reported in this paper.

%% ============================================================
%%  Data availability
%% ============================================================
\section*{Code availability}
The code used to produce the numerical results in this paper is
implemented in the open-source Firedrake finite element framework
(\url{https://www.firedrakeproject.org}) and is available at
\url{https://github.com/karinakowalczyk/swe_projection}, including
all scripts required to reproduce the test cases.

%% ============================================================
%%  Acknowledgements
%% ============================================================
%\section*{Acknowledgements}
% Fill in funding sources, grants, and personal acknowledgements.

%% ============================================================
%%  Bibliography
%% ============================================================
\bibliographystyle{elsarticle-num}
\bibliography{references}

@phdthesis{hastermann2022analysis,
  title={Analysis of the Cell-Vertex Finite Volume Method for Pseudo-Incompressible Divergence Constraints on Quadrilateral and Cuboid Meshes},
  author={Hastermann, Gottfried},
  year={2022}
}

@article{bernard2009high,
  title={High-order discontinuous Galerkin schemes on general 2D manifolds applied to the shallow water equations},
  author={Bernard, P-E and Remacle, J-F and Comblen, Richard and Legat, Vincent and Hillewaert, Koen},
  journal={Journal of Computational Physics},
  volume={228},
  number={17},
  pages={6514--6535},
  year={2009},
  publisher={Elsevier}
}

@article{arnold2006finite,                                                      
  title={Finite element exterior calculus, homological techniques, and applicat\
ions},                                                                          
  author={Arnold, Douglas N and Falk, Richard S and Winther, Ragnar},           
  journal={Acta Numerica},                                                      
  volume={15},                                                                  
  pages={1--155},                                                               
  year={2006},                                                                  
  publisher={Cambridge University Press}                                        
}

@book{arnold2018finite,                                                         
  title={Finite element exterior calculus},                                     
  author={Arnold, Douglas N},                                                   
  year={2018},                                                                  
  publisher={SIAM}                                                              
}

@article{arnold2010finite,                                                      
  title={Finite element exterior calculus: from {H}odge theory to Numerical     
  stability},                                                                   
  author={Arnold, Douglas and Falk, Richard and Winther, Ragnar},               
  journal={Bulletin of the American Mathematical Society},                      
  volume={47},                                                                  
  number={2},                                                                   
  pages={281--354},                                                             
  year={2010}                                                                   
}

@article{robert1972implicit,
  author    = {A. Robert and J. Henderson and C. Turnbull},
  title     = {An Implicit Time Integration Scheme for Baroclinic Models of the Atmosphere},
  journal   = {Monthly Weather Review},
  volume    = {100},
  number    = {5},
  pages     = {329--335},
  year      = {1972}
}

@article{klemp1978numerical,
  author  = {Klemp, J. B. and Wilhelmson, R. B.},
  title   = {The simulation of three-dimensional convective
             storm dynamics},
  journal = {Journal of the Atmospheric Sciences},
  volume  = {35},
  number  = {6},
  pages   = {1070--1096},
  year    = {1978}
}

@article{skamarock1992efficiency,
  author  = {Skamarock, W. C. and Klemp, J. B.},
  title   = {The Stability of Time-Split Numerical Methods for the Hydrostatic and the Nonhydrostatic Elastic Equations},
  journal = {Monthly Weather Review},
  year    = {1992},
  volume  = {120},
  number  = {9},
  pages   = {2109--2127}
}

@article{weller2013runge,
  author  = {Weller, H. and Lock, S.-J. and Wood, N.},
  title   = {{Runge--Kutta IMEX} schemes for the horizontally
             explicit/vertically implicit ({HEVI}) solution of
             wave equations},
  journal = {Journal of Computational Physics},
  volume  = {252},
  pages   = {365--381},
  year    = {2013}
}

@article{lock2014numerical,
  author  = {Lock, S.-J. and Wood, N. and Weller, H.},
  title   = {Numerical analyses of {Runge--Kutta} implicit--explicit
             schemes for horizontally explicit, vertically implicit
             solutions of atmospheric models},
  journal = {Quarterly Journal of the Royal Meteorological Society},
  volume  = {140},
  number  = {682},
  pages   = {1654--1669},
  year    = {2014}
}

@article{vater2009stability,
  author  = {Vater, S. and Klein, R.},
  title   = {Stability of a {Cartesian} grid projection method
             for zero {Froude} number shallow water flows},
  journal = {Numerische Mathematik},
  volume  = {113},
  number  = {1},
  pages   = {123--161},
  year    = {2009}
}

@article{vater2018semi,
  author  = {Vater, S. and Beisiegel, N. and Behrens, J.},
  title   = {A semi-implicit multiscale scheme for shallow water
             flows at low {Froude} number},
  journal = {Communications in Applied Mathematics and
             Computational Science},
  volume  = {13},
  number  = {2},
  pages   = {303--336},
  year    = {2018}
}

@article{benacchio2019semi,
  author  = {Benacchio, T. and Klein, R.},
  title   = {A semi-implicit compressible model for atmospheric
             flows with seamless access to soundproof and
             hydrostatic dynamics},
  journal = {Monthly Weather Review},
  volume  = {147},
  number  = {11},
  pages   = {4221--4240},
  year    = {2019}
}

@article{benacchio2014blended,
  author  = {Benacchio, T. and O'Neill, W. P. and Klein, R.},
  title   = {A blended soundproof-to-compressible numerical
             model for small- to mesoscale atmospheric dynamics},
  journal = {Monthly Weather Review},
  volume  = {142},
  number  = {12},
  pages   = {4416--4438},
  year    = {2014}
}

@article{klein2016doubly,
  author  = {Klein, R. and Benacchio, T.},
  title   = {A doubly blended model for multiscale atmospheric
             dynamics},
  journal = {Journal of the Atmospheric Sciences},
  volume  = {73},
  number  = {3},
  pages   = {1179--1186},
  year    = {2016}
}

@article{chew2022one,
  author  = {Chew, R. and Benacchio, T. and Hastermann, G.
             and Klein, R.},
  title   = {A one-step blended soundproof-compressible model
             with balanced data assimilation},
  journal = {Monthly Weather Review},
  volume  = {150},
  number  = {9},
  pages   = {2231--2254},
  year    = {2022}
}

@article{cotter2012mixed,
  author  = {Cotter, C. J. and Shipton, J.},
  title   = {Mixed finite elements for numerical weather
             prediction},
  journal = {Journal of Computational Physics},
  volume  = {231},
  number  = {21},
  pages   = {7076--7091},
  year    = {2012}
}

@article{cotter2014finite,
  author  = {Cotter, C. J. and Thuburn, J.},
  title   = {A finite element exterior calculus framework for
             the rotating shallow-water equations},
  journal = {Journal of Computational Physics},
  volume  = {257},
  pages   = {1506--1526},
  year    = {2014}
}

@article{natale2016compatible,
  title={Compatible finite element spaces for geophysical fluid dynamics},
  author={Natale, Andrea and Shipton, Jemma and Cotter, Colin J},
  journal={Dynamics and Statistics of the Climate System},
  volume={1},
  number={1},
  pages={dzw005},
  year={2016},
  publisher={Oxford University Press UK}
}

@article{melvin2019mixed,
  title={A mixed finite-element, finite-volume, semi-implicit discretization for atmospheric dynamics: Cartesian geometry},
  author={Melvin, Thomas and Benacchio, Tommaso and Shipway, Ben and Wood, Nigel and Thuburn, John and Cotter, Colin},
  journal={Quarterly Journal of the Royal Meteorological Society},
  volume={145},
  number={724},
  pages={2835--2853},
  year={2019},
  publisher={Wiley Online Library}
}

@article{melvin2024mixed,
  title={A mixed finite-element, finite-volume, semi-implicit discretisation for atmospheric dynamics: Spherical geometry},
  author={Melvin, Thomas and Shipway, Ben and Wood, Nigel and Benacchio, Tommaso and Bendall, Thomas and Boutle, Ian and Brown, Alex and Johnson, Christine and Kent, James and Pring, Stephen and others},
  journal={Quarterly Journal of the Royal Meteorological Society},
  volume={150},
  number={764},
  pages={4252--4269},
  year={2024},
  publisher={Wiley Online Library}
}

@article{melvin2010inherently,
  title={An inherently mass-conserving iterative semi-implicit semi-Lagrangian discretization of the non-hydrostatic vertical-slice equations},
  author={Melvin, Thomas and Dubal, Mark and Wood, Nigel and Staniforth, Andrew and Zerroukat, Mohamed},
  journal={Quarterly Journal of the Royal Meteorological Society},
  volume={136},
  number={648},
  pages={799--814},
  year={2010},
  publisher={Wiley Online Library}
}

@misc{gusto2023,
  author       = {{Gusto Development Team}},
  title        = {Gusto: a compatible finite element dynamical
                  core toolkit},
  howpublished = {\url{https://www.firedrakeproject.org/gusto}},
  year         = {2023}
}

@article{cockburn2001runge,
  author  = {Cockburn, B. and Shu, C.-W.},
  title   = {{Runge--Kutta} discontinuous {Galerkin} methods
             for convection-dominated problems},
  journal = {Journal of Scientific Computing},
  volume  = {16},
  number  = {3},
  pages   = {173--261},
  year    = {2001}
}

@article{ern2022invariant,
  author  = {Ern, A. and Guermond, J.-L.},
  title   = {Invariant-domain preserving high-order time
             stepping: {I.} Explicit {Runge--Kutta} schemes},
  journal = {SIAM Journal on Scientific Computing},
  volume  = {44},
  number  = {5},
  pages   = {A3366--A3392},
  year    = {2022}
}

@article{cockburn2009unified,
  author  = {Cockburn, B. and Gopalakrishnan, J. and Lazarov, R.},
  title   = {Unified hybridization of discontinuous {Galerkin},
             mixed, and continuous {Galerkin} methods for second
             order elliptic problems},
  journal = {SIAM Journal on Numerical Analysis},
  volume  = {47},
  number  = {2},
  pages   = {1319--1365},
  year    = {2009}
}

@article{brezzi1985two,
  author  = {Brezzi, F. and Douglas, J. and Marini, L. D.},
  title   = {Two families of mixed finite elements for second
             order elliptic problems},
  journal = {Numerische Mathematik},
  volume  = {47},
  number  = {2},
  pages   = {217--235},
  year    = {1985}
}

@manual{FiredrakeUserManual,
  title        = {Firedrake User Manual},
  author       = {David A. Ham and Paul H. J. Kelly and Lawrence Mitchell and Colin J. Cotter and Robert C. Kirby and Koki Sagiyama and Nacime Bouziani and Sophia Vorderwuelbecke and Thomas J. Gregory and Jack Betteridge and Daniel R. Shapero and Reuben W. Nixon-Hill and Connor J. Ward and Patrick E. Farrell and Pablo D. Brubeck and India Marsden and Thomas H. Gibson and Miklós Homolya and Tianjiao Sun and Andrew T. T. McRae and Fabio Luporini and Alastair Gregory and Michael Lange and Simon W. Funke and Florian Rathgeber and Gheorghe-Teodor Bercea and Graham R. Markall},
  organization = {Imperial College London and University of Oxford and Baylor University and University of Washington},
  edition      = {First edition},
  year         = {2023},
  month        = {5},
  doi          = {10.25561/104839},
}

@article{williamson1992standard,
  author  = {Williamson, D. L. and Drake, J. B. and Hack, J. J.
             and Jakob, R. and Swarztrauber, P. N.},
  title   = {A standard test set for numerical approximations
             to the shallow water equations in spherical geometry},
  journal = {Journal of Computational Physics},
  volume  = {102},
  number  = {1},
  pages   = {211--224},
  year    = {1992}
}

@article{galewsky2004initial,
  author  = {Galewsky, J. and Scott, R. K. and Polvani, L. M.},
  title   = {An initial-value problem for testing numerical
             models of the global shallow-water equations},
  journal = {Tellus~A},
  volume  = {56},
  number  = {5},
  pages   = {429--440},
  year    = {2004}
}

@article{lauter2005unsteady,
  author  = {L{\"a}uter, M. and Handorf, D. and Dethloff, K.},
  title   = {Unsteady analytical solutions of the spherical
             shallow water equations},
  journal = {Journal of Computational Physics},
  volume  = {210},
  number  = {2},
  pages   = {535--553},
  year    = {2005}
}

@book{hairer2006geometric,
  author    = {Hairer, E. and Lubich, C. and Wanner, G.},
  title     = {Geometric Numerical Integration:
               Structure-Preserving Algorithms for Ordinary
               Differential Equations},
  edition   = {2nd},
  publisher = {Springer},
  address   = {Berlin},
  year      = {2006}
}

@book{durran2010numerical,
  author    = {Durran, D. R.},
  title     = {Numerical Methods for Fluid Dynamics:
               With Applications to Geophysics},
  edition   = {2nd},
  publisher = {Springer},
  address   = {New York},
  year      = {2010}
}

\end{document}